\documentclass[a4paper, english, 12pt, reqno, dvipsnames]{amsart}

\usepackage[english]{babel}
\usepackage[utf8]{inputenc}
\usepackage{a4wide}

\usepackage{graphics}
\graphicspath{ {./images/} }

\usepackage{tabularx,booktabs}
\usepackage[table]{xcolor}

\usepackage{bm}
\usepackage{pgfplots}
\pgfplotsset{compat=1.18}
\usepackage{algorithm,algorithmic}
\usepackage{dsfont}
\usepackage{diagbox}

\usepackage[normalem]{ulem}
\newcounter{AR}
\definecolor{violet}{rgb}{0.580,0.,0.827}
\newcommand{\AR}[3]{\typeout{Warning : a correction remains in page \thepage}\stepcounter{AR}%
{\color{blue}\ifmmode\text{\,\sout{\ensuremath{#1}}\,}\else\sout{#1}\fi}%
{\color{red}#2}{ [\color{violet} AR: #3]}}

\newcommand{\displ}{\mathbf{u}} %displacement
\newcommand{\stress}{\bm{\sigma}} %stress tensor
\newcommand{\strain}{\bm{\varepsilon}}
\newcommand{\strainInd}{\varepsilon}
\newcommand{\x}{ \mathbf{x}} 
 
\newcommand{\R}{\mathbb{R}}

\newcommand{\norm}[2]{\| #1 \|_{#2}}
\newcommand{\code}[1]{%
  \begingroup
  \ttfamily
  \begingroup\lccode`~=`/\lowercase{\endgroup\def~}{/\discretionary{}{}{}}%
  \begingroup\lccode`~=`[\lowercase{\endgroup\def~}{[\discretionary{}{}{}}%
  \begingroup\lccode`~=`.\lowercase{\endgroup\def~}{.\discretionary{}{}{}}%
  \begingroup\lccode`~=`_\lowercase{\endgroup\def~}{_\discretionary{}{}{}}%
  \catcode`/=\active\catcode`[=\active\catcode`.=\active\catcode`_=\active
  \scantokens{#1\noexpand}%
  \endgroup
}
\numberwithin{equation}{section}

\makeatletter \@addtoreset{table}{section} \makeatother

\usepackage[colorlinks = true, linkcolor = blue, citecolor = blue, urlcolor = blue]{hyperref}

\title[Recovery of elastic material parameters in induction motor rotors]{Recovery of transversely-isotropic elastic material parameters in induction motor rotors}

\author[H.M.~Cheng]{Hanz M. Cheng}
\address{School of Engineering Science, Lappeenranta--Lahti University of Technology, P.O. Box 20, 53851 Lappeenranta, Finland}
\email{hanz.cheng@lut.fi}

\author[T.~Helin]{Tapio Helin}
\address{School of Engineering Science, Lappeenranta--Lahti University of Technology, P.O. Box 20, 53851 Lappeenranta, Finland}
\email{tapio.helin@lut.fi}

\author[V.-P.~Manninen]{Ville-Petteri Manninen}
\address{School of Engineering Science, Lappeenranta--Lahti University of Technology, P.O. Box 20, 53851 Lappeenranta, Finland}
\email{ville-petteri.manninen@student.lut.fi}
 
\author[T.~Holopainen]{Timo Holopainen}
\address{Large Motors and Generators, ABB Oy, P.O. Box 186, FI-00381 Helsinki, Finland}
\email{timo.holopainen@fi.abb.com}

\author[J.~Jokinen]{Juha Jokinen}
\address{Large Motors and Generators, ABB Oy, P.O. Box 186, FI-00381 Helsinki, Finland}
\email{juha.jokinen@fi.abb.com}

\author[S.~Sorvari]{Samu Sorvari}
\address{Large Motors and Generators, ABB Oy, P.O. Box 186, FI-00381 Helsinki, Finland}
\email{samu.sorvari1@fi.abb.com}
	
\author[A.~Rupp]{Andreas Rupp}
\address{School of Engineering Science, Lappeenranta--Lahti University of Technology, P.O. Box 20, 53851 Lappeenranta, Finland}
\email{andreas.rupp@fau.de}
	
\thanks{%
This work has been supported by the Research Council of Finland (RCoF) through the \emph{Flagship of advanced mathematics for sensing, imaging and modelling}, decision number 358944. Moreover, TH and AR were supported though RCoF decision numbers 353094, 348504, 326961, 350101, 354489 and 359633. In addition, AR was supported by Business Finland's project number 539/31/2023.%
}

\begin{document}

\begin{abstract}
 We propose numerical algorithms for recovering parameters in eigenvalue problems for linear elasticity of transversely isotropic materials. Specifically, the algorithms are used to recover the elastic constants of a rotor core. Numerical tests show that in the noiseless setup, two pairs of bending modes are sufficient for recovering one to four parameters accurately. To recover all five parameters that govern the elastic properties of electric engines accurately, we require three pairs of bending modes and one torsional mode. Moreover, we study the stability of the inversion method against multiplicative noise; for tests in which the data contained multiplicative noise of at most 1\%, we find that all parameters can be recovered with an error less than $10\%$.  
 \\[1ex] \noindent \textsc{Keywords.}
 Ensemble Kalman inversion, inverse eigenvalue problem, least squares minimization, linear elasticity, modal testing, and transversely isotropic material. 
\end{abstract}

\date{\today}

\maketitle

% ---------------------------------------------------------------------------
\section{Introduction}
% ---------------------------------------------------------------------------

% 
%\tapio{(Paragraph on motivation)}
The rotor plays a vital role in an electric motor's operation and efficiency. Due to increasing computational resources, accurate mathematical modeling of the rotor's properties is becoming more feasible with finite element analysis and multi-physics simulations. However, due to the rotor's complex structure, the manufacturing process will alter the rotor's elastic parameters to some effect, influencing, e.g., performance, resonant frequencies, and durability of the motor. Therefore, it is helpful for the manufacturer to verify that the produced rotor meets the design specifications experimentally. In this paper, we develop novel numerical methods for this purpose.

%\tapio{(Explain modal testing)}
The standard experimental technique to study the rotor properties is modal testing \cite{E00-testing}, which is used to determine dynamic properties such as natural frequencies, mode shapes, and damping ratios. It is performed by applying a controlled excitation to the rotor core and measuring the corresponding response. The excitation can be delivered in various forms, including impact hammers or electromagnetic shakers, to induce vibrations. Sensors, typically accelerometers, are placed on the structure to measure the vibrations and capture the dynamic response across a range of frequencies. Following this, the measured data is typically analyzed by subjecting it to a range of curve-fitting procedures to find the mathematical model that provides the closest description of the observed measurements. 

\begin{figure}[t!]\centering
 \includegraphics[width=0.75\linewidth]{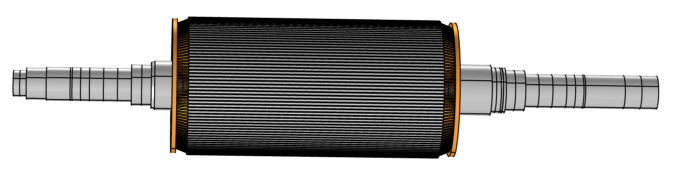}
 \begin{tabular}{ccc}
  \raisebox{1\height}{\includegraphics[scale=0.25]{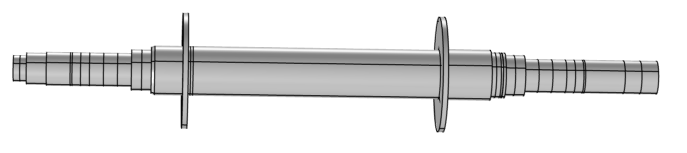}} &
  \raisebox{.075\height}{\includegraphics[width=0.30\linewidth]{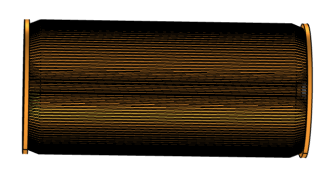}} &
  \includegraphics[width=0.30\linewidth]{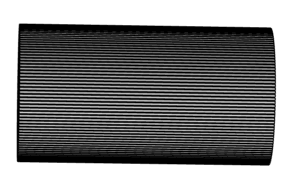}
 \end{tabular}
 \caption{Sketch of an induction motor's rotor. The top picture shows an assembled rotor consisting of a shaft with pressure plates, short circuit rings with bars, and a rotor core, which are depicted individually from left to right in the bottom row.}\label{fig.rotor_components}
\end{figure}

%\tapio{(Explain what rotor components and what are the parameters of interest)}
In this paper, we utilize the spectral data obtained in modal testing for the identification of the stiffness tensor of the rotor core by optimization. To be precise, we consider a rotor composed of three different types of materials (see Figure \ref{fig.rotor_components}). Firstly, the shaft and pressure plates are modeled using the structural steel isotropic material model. Secondly, the short circuit rings and the bars are modeled using the isotropic material model of copper. The material properties of these two parts are assumed to be known. Thirdly, the rotor core is a stack of thin electrical steel sheets pressed together by pressure plates on both sides of the core. In addition, the steel sheets are fixed by swaging longitudinally the rotor bars between the core teeth, forming transversely isotropic material properties specified by the symmetry axis along the shaft. Here, we assume that the stiffness tensor is
% \tapio{(I dropped "almost" for clarity)}
constant within the rotor core, and the parameters describing it are partially or completely unknown and, therefore, are our parameters of interest.

% \tapio{We should address in more detail here, what is the added value of identifying the stiffness tensor (instead of just the eigenfrequencies)}
Ultimately, our approach allows engineers to predict the eigenfrequencies of similar but different rotors, which they design virtually. These eigenfrequencies are critically important for rotors as they should not be run with a rotation speed that resembles their eigenfrequencies, which might damage the electric engine. Thus, engineers can design useful engines virtually, without building the engines and checking their eigenfrequencies, which significantly reduces the development costs.

We employ linear elasticity to characterize the dynamic response of the rotor, leading to an inverse eigenvalue problem of the associated elliptic partial differential equation 
\cite{kachalov2001inverse}. This equation encompasses a stiffness tensor that describes the material's physical parameters. For the transversely isotropic material in the rotor core, the stiffness tensor is specified by five parameters composed of Young's modulus, shear modulus, and the Poisson ratio (the first two described by two parameters each; one for the symmetry plane and one for the symmetry axis). Each of the five parameters is assumed to be constant throughout the rotor core. Assuming a multiplicative noise model in the measurement data, we then treat the recovery of these parameters as an optimization problem of the maximum likelihood density.

The inverse problems literature considering parameter identification problems in linear elasticity typically aims at recovering a spatially varying stiffness tensor based on Cauchy data on the boundary or some scattering data \cite{kachalov2001inverse}. The theoretical problem is notoriously tricky and open in the general case (see relevant literature below). Moreover, inverse eigenvalue problems such as here with only spectral data available, tend not to possess unique solutions even in most superficial elliptic problems. Instead, multiple spectra corresponding to varying boundary conditions are usually observed to guarantee uniqueness \cite{kachalov2001inverse}. This paper considers a substantially simpler setup, where the stiffness tensor is constant in the rotor core, begging the question of whether partial recovery of the parameters is possible.

% ---------------------------------------------------------------------------
\subsection{Our contribution}
% ---------------------------------------------------------------------------
% 
This paper studies two optimization approaches for the inverse eigenvalue problem: state-of-the-art Hessian-based methods called Sequential Least-Squares Programming (SLSQP) and a derivative-free method called ensemble Kalman inversion (EKI). Our numerical studies are carried out exclusively with simulated data, and the application to real-life data is left part of future studies. However, we limit our experiments to the realistic situation where only 4 to 7 lowest eigenfrequencies are observed. In particular, we utilize EKI to conclude the identifiability of parameters. Our observations are as follows:
\begin{itemize}
 \item The Hessian-based methods are relatively fast in terms of computational cost, but identifying the parameters to a good approximation is challenging. The second approach with EKI requires a substantially higher computational cost but converges more reliably. This is demonstrated in Sections \ref{sec:lsq_nonoise} -- \ref{sec:EKI_nonoise}.
 \item The EKI method identifies the stiffness tensor reasonably well from noise-free data. Moreover, we compare the identifiability of each single parameter and parameter pairs (see Tables \ref{tab:noNoise1} and \ref{tab:noNoise2}) with SLSQP from different numbers of measured eigenvalues. The results indicate that EKI generally approximates the parameters better than the SLSQP. We also note that the Poisson ratio is most sensitive to numerical errors as it has a negligible impact on the eigenvalues, as noted under the related setting in \cite{DQetal23-identification}.
 \item We study the stability of the inverse problem with the EKI method by comparing reconstruction error with decreasing noise level. We consider subproblems where sets of two, three, and five parameters are identified. Moreover, we study the influence of the number of measured eigenfrequencies. We find that the entirely inverse problem is severely ill-posed; e.g., identifying the four main parameters (excluding Poisson ratio) up to the relative error of order $10^{-2}$ requires the factor of the multiplicative noise to be less than $10^{-5}$. Moreover, recovering the Poisson ratio with the same tolerance requires nearly perfect knowledge of the spectrum.
\end{itemize}

% ---------------------------------------------------------------------------
\subsection{Literature review}
% ---------------------------------------------------------------------------

In the mathematical literature on inverse problems, the main uniqueness results with Cauchy data and spatially dependent stiffness tensors were developed in \cite{nakamura1994global, nakamura2003global} and \cite{eskin2002inverse, eskin2004inverse}, where the authors considered a fully isotropic material and demonstrate unique recovery of the spatially dependent Lam\'{e} parameters. For the various aspects of inverse problems with elastic medium, see \cite{bonnet2005inverse}. Let us note that there is also recent work on anisotropic medium \cite{de2023reconstruction}, but the general case is poorly understood. 

Eigenvalue problems appear in numerous applications of computational science and engineering. Here, let us mention vibration analysis of mechanical systems (see \cite{mottershead2006inverse} and references therein). Accordingly, a large body of mathematical literature has developed around eigenvalue perturbation analysis, dating back to seminal works by Rellich and Kato \cite{rellich1969perturbation, kato2013perturbation}. There is a broad literature on inverse eigenvalue problems; see, e.g., \cite{yamamoto1990inverse, rundell1992reconstruction, chu1998inverse, horvath2005inverse}. For a general overview, see \cite{kachalov2001inverse, chadan1997introduction}.

The rotor of an induction motor is a fabricated structure consisting of hundreds of steel laminates and other parts. These parts are fitted together by different connection methods like friction joints, interference fits, swaging, welding, and soldering. In addition, some of the fabricated rotors are resin treated, which affects the joints of the structure. The connection of parts produces a solid and stable entity against mechanical loads like centrifugal force. Finally, due to the different thermal expansion factors of steel and winding materials like copper, the connections between the parts may change as function of rotor temperature. The accurate modelling of induction rotors, starting from separate parts and their connections, is out of the reach of current state-of-the-art. Thus, various methods have been developed to simplify the rotor core models and to identify the corresponding equivalent material properties. These methods have been based on the modal testing of free-free rotors \cite{Kim06-lamination,Seo10,Singhal2011EffectOL}.
    
The main objective of this work is to recover, via modal testing, a mathematical model that provides the closest description of the observed measurements, so that it can be used, e.g. for engineers to design engines virtually, without having to build the engines and checking their eigenfrequencies. In this paper, we consider a linear elasticity model, for which the recovery of the mathematical model translates to recovering the material properties associated with the stiffness tensor. The recovery of these material properties may involve direct methods stemming from homogenization, such as those presented in \cite{MillithalerEtal15-homogenization,BZFW20-identify} and the references therein. Closer to this work are approaches which fits the response of the numerical model to the corresponding experimental result in the least squares sense by minimizing an error estimator \cite{LBAK16-dynamic_approach}, illustrated for example via the Nelder-Mead method \cite[Chapter 4.2.3]{R09-thesis}, the feasible directions non-linear interior point algorithm \cite{AraujoEtal00-numerical_experimental}, Newton-type methods \cite{SaitoEtal16-ortho_identification}, or CNN (convolutional neural network)-MPGA (multi-population genetic algorithm) \cite{DQetal23-identification}.

% seen for example in \cite{AraujoEtal00-numerical_experimental,DQetal23-identification,LBAK16-dynamic_approach,R09-thesis,SaitoEtal16-ortho_identification}

With regards to our optimization algorithms, we note that the Hessian-based methods are standard, and a detailed exposition is provided, e.g. in \cite{NW06}. Here, the numerical estimation of the Hessian is challenging, and derivative-free methods provide a tempting option. Ensemble Kalman methods were initially introduced for dynamic state estimation in \cite{evensen2009data}. In \cite{chen2012ensemble, emerick2013investigation}, these methods also serve as derivative-free solvers for Bayesian inverse problems. Further modification in \cite{I13-EKI_IP} transformed these methods into derivative-free optimizers, which are often referred to as the Ensemble Kalman Inversion (EKI). While noting that many other general-purpose derivative-free optimization methods exist, the EKI method was chosen here as it has been shown, both theoretically in simple settings and empirically, to perform very well in the context of inverse problems \cite{I13-EKI_IP, schillings2017analysis, schillings2018convergence, ding2021ensemble} and has generated a number of variants.

% ---------------------------------------------------------------------------
\subsection{Structure of this work} 
% ---------------------------------------------------------------------------
% 
 Section \ref{sec:setup_and_model} describes the experimental setup for obtaining the data, together with the eigenvalue problem and the mathematical model, which links the measurement data to the parameters of interest. Afterward, Section \ref{sec:num_methods} illuminates the inverse problem and some algorithms for solving the inverse problem. In particular, we explore a least squares minimization approach and ensemble Kalman inversion (EKI) \cite{I13-EKI_IP,I16_regularising_EKI} for recovering the parameters.  Numerical tests and results are presented in Section \ref{sec:Num_res}. A summary of the numerical results, the methods' usage and limitations, and recommendations for future work can be found in Section \ref{sec:summary}.
\section{Measurements and mathematical model} \label{sec:setup_and_model}
% ---------------------------------------------------------------------------
% 
In a typical experiment, the rotor is laid at rest. Sensors are then strategically placed on the rotor at a few (typically three to five) locations. Subsequently, the rotor is hit with a hammer, and the resulting vibrations are captured by a measurement device attached to the sensors. From these vibrations, the eigenfrequencies are extracted. Afterward, the measured data is analyzed, and we select an appropriate mathematical model that relates the parameters to the eigenfrequencies. 

The measurement data obtained from the experiment correspond to the spectral values, i.e., eigenvalues and eigenmodes of a vibrating engine. Thus, we want to use a mathematical model that relates the same spectral values to the elastic parameters. We say the parameter values are good candidates for accurate physical material parameters if they yield the measured eigenvalues and eigenmodes if used in the mathematical model.

% ---------------------------------------------------------------------------
\subsection{The partial differential equations}
Let us denote the spatial domain of the engine at rest by $\Omega \subset \mathds R^3$. The domain $\Omega = \Omega_\mathrm{c} \cup \Omega_\mathrm{s} \cup \Omega_\mathrm{r}$ comprises three disjoint components: the copper component $\Omega_\mathrm{c}$, the steel component $\Omega_\mathrm{s}$, and the component consisting of transversely isotropic material $\Omega_\mathrm{r}$. Function $\displ$ models the elastic displacement the rotor experiences when excitation is delivered. We model the relation of a vibrating rotor to its elastic constants as the eigenvalue problem for linear elasticity. 
By Newton's second law, the stress tensor $\stress$ satisfies, in the absence of force, that
% Newton's second law on a domain $\Omega$ over the time interval $(0,T)$ can be written as 
% \begin{equation}\label{eq:Newton_stress}
% 	\nabla \cdot \stress +\force = \rho\frac{\partial^2 \displ}{\partial t^2} \quad \mathrm{on} \quad  \Omega\times(0,T).
% \end{equation}
\begin{equation}\label{eq:Newton_stress}
 \nabla \cdot \stress = \rho\frac{\partial^2 \displ}{\partial t^2} \quad \mathrm{in} \quad  \Omega\times(0,T).
\end{equation}
where $\rho$ denotes the mass density and $\displ$ describes the displacement vector. For this paper, we consider a linear elasticity model, where the stress is related to the strain via Hooke's law
\begin{equation} \label{eq:Hookes}
\begin{aligned}
    \stress = \begin{cases}
        \mathbf{C}_c \strain \quad \mathrm{in} \quad \Omega_c,\\
        \mathbf{C}_s \strain \quad \mathrm{in} \quad \Omega_s, \\
        \mathbf{C}_\Gamma \strain \quad \mathrm{in} \quad \Omega_\Gamma
    \end{cases}
\end{aligned}
\end{equation}
for some fourth-order stiffness tensors $\mathbf{C}_r$, $r\in\{c,s,\Gamma\}$ and strain tensor $\strain$. The kinematic equation governs the strain tensor $\strain$
\begin{equation} \label{eq:strain}
	\strain = \frac{1}{2}\big(\nabla \displ + \nabla \displ^\top \big).
\end{equation}
For time-harmonic solutions of the form $\displ(\x,t) = \displ(\x)e^{i\omega t}$, we obtain the eigenvalue problem
\begin{align}
   \omega^2 \rho  \displ + \nabla \cdot \stress &= 0 \quad \mathrm{in} \quad \Omega,\label{eq:eigenvalue_problem}
\end{align}
where the eigenvalues are denoted by $\Lambda = \omega^2$.
% with boundary conditions
% \begin{equation} \label{eq:BCs}
% \displ = 0 \quad \mathrm{on} \,\,\, \partial\Omega_D, \qquad \stress \mathbf{n}=0 \quad \mathrm{on} \,\,\, \partial\Omega_N,
% \end{equation}
% where we assume that the boundary is composed of two disjoint components $\partial \Omega_D$ and $\partial \Omega_N$, i.e., $\partial\Omega = \partial \Omega_D \bigcup \partial \Omega_N$. 

In this paper, we consider unconstrained structures, i.e. the elastic body is not subject to any geometric constraints. Since we are studying a three-dimensional model, there are six rigid body modes (corresponding to zero eigenvalues) \cite{E00-testing}. It is also well-known that the weak form of the eigenvalue problem is well-posed, and therefore we can use, with justification, various numerical methods based on the weak form, such as the FEM \cite{EL13_existence}. We also note that the mathematical model is mainly used for performing an eigenfrequency analysis, i.e., for recovering the shape of the eigenmodes. The amplitude of the physical vibrations can only be determined if an actual excitation is known together with the damping properties, if any. As such, one may choose eigenmodes of unit amplitude to obtain unique solutions. 

% ---------------------------------------------------------------------------
\subsection{Stress-strain relations}
% ---------------------------------------------------------------------------
% 
This section outlines an explicit relationship between the stress $\stress$ and the strain $\strain$, employing conventional notation. Since the strain tensor 
%% personal note: this is different from the std engineering notation in the sense that 4 and 6 have been switched, so upon comparison with literature, the entries in the 4th and 6th rows of the matrix are interchanged
\begin{equation}
    \strain = \begin{bmatrix}
    \strainInd_{11} & \strainInd_{12} & \strainInd_{13} \\
    \strainInd_{21} & \strainInd_{22} & \strainInd_{23} \\
    \strainInd_{31} & \strainInd_{32} & \strainInd_{33} \\
    \end{bmatrix}
\end{equation}
is symmetric, we denote 
\begin{equation} \nonumber
    \begin{aligned}
        \bar\strainInd_1 &= \strainInd_{11},  \quad \bar\strainInd_2 &= \strainInd_{22},  \quad \bar\strainInd_3 &= \strainInd_{33}, \\
        \bar\strainInd_4 &= \strainInd_{23},  \quad \bar\strainInd_5 &= \strainInd_{13}, 
        \quad \bar\strainInd_6 &= \strainInd_{12}. \\
    \end{aligned}
\end{equation}
By using a similar notation for the stress tensor,
% Similarly, for the stress tensor $\stress$, we denote
% \begin{equation} \nonumber
%     \begin{aligned}
%         \bar\stress_1 &= \stress_{11},  \quad \bar\stress_2 &= \stress_{22},  \quad \bar\stress_3 &= \stress_{33}, \\
%         \bar\stress_4 &= \stress_{12},  \quad \bar\stress_5 &= \stress_{13},  \quad \bar\stress_6 &= \stress_{23}. \\
%     \end{aligned}
% \end{equation}
we now have $\bar{\stress}, \bar{\strain} \in \R^{6\times 1}$ and rewrite Hooke's law as 
\begin{equation}\label{eq:Hookes_stiff}
\bar\stress = \mathbf{\bar{C}} \bar\strain.
\end{equation}
Here, the stress is related to the strain via the stiffness tensor $\mathbf{\bar{C}}$, which is why \eqref{eq:Hookes_stiff} is referred to as Hooke's law in \emph{stiffness form}. Equivalently, we may write Hooke's law in \emph{compliance form}
\begin{equation}\label{eq:Hookes_comp}
\bar\strain = \mathbf{\bar{S}} \bar\stress,
\end{equation}
where $\mathbf{\bar{S}} = \mathbf{\bar{C}}^{-1}$ is known as the compliance tensor. As seen in equation \eqref{eq:Hookes}, the tensor depends on the type of material being considered. We present some of the different forms taken by the compliance tensor $\mathbf{\bar{S}}$ for isotropic and transversely isotropic materials. For a more detailed discussion, we refer the readers to \cite{P04-formula-matrices,S21-lin_el}. 

% ---------------------------------------------------------------------------
\subsubsection{Isotropic material}
% ---------------------------------------------------------------------------
% 
Isotropic materials exhibit uniform material properties regardless of direction. In other words, their compliance tensor is symmetric and depends on only two parameters: Young's modulus $E$ and Poisson ratio $\nu$. It can be formulated as
\begin{equation}
    \mathbf{\bar{S}} = \frac{1}{E}\begin{bmatrix}
    1 & -\nu & -\nu & 0 & 0 & 0 \\
    -\nu & 1 & -\nu & 0 & 0 & 0 \\
    -\nu & -\nu & 1 & 0 & 0 & 0 \\
    0 & 0 & 0 & 1+\nu & 0 & 0 \\
    0 & 0 & 0 & 0 & 1+\nu & 0 \\
     0 & 0 & 0 & 0 & 0 & 1+\nu
    \end{bmatrix}.
\end{equation}
For a more detailed analysis of the eigenvalue problem, obtaining an explicit form that relates stress and strain can be useful. For isotropic materials, we may write this relation as
\begin{equation} \label{eq:stress_strain_isotropic}
	\stress = 2G\strain + \lambda \mathrm{tr}(\strain) \mathbf{I},
\end{equation}
where $G,\lambda$ are Lam\'e constants, $\mathrm{tr}(\strain)$ is the trace of the strain tensor and $\mathbf{I}$ is the identity matrix. The Lam\'e constant $G$ pertains to the shear modulus, which can be obtained via the relation \cite{LLS60-shear}
\begin{equation}\nonumber
E = 2G(1+\nu),
\end{equation}
and the Lam\'e constant $\lambda$ is related to the Young's modulus and Poisson's ratio by the following equation:
\begin{equation}
	\begin{aligned} 
		\lambda &= \frac{E}{1-2\nu}\frac{\nu}{1+\nu}.
	\end{aligned} 
\end{equation}
We note, however, that the expression \eqref{eq:stress_strain_isotropic} is not well-defined when $\lambda \rightarrow \infty$, i.e., when the material is incompressible \cite{CS04-FE}. Hence, it is more useful to write the compliance form 
\begin{equation} \label{eq:strain_stress_isotropic}
	\strain = \frac{1}{2G}\left(\stress-\frac{\lambda}{3\lambda+2G}\mathrm{tr}(\stress)\mathbf{I}\right).
\end{equation}

% ---------------------------------------------------------------------------
\subsubsection{Transversely isotropic material}
% ---------------------------------------------------------------------------
% 
Transversely isotropic materials have the same properties in one plane, called the plane of isotropy and different properties in the direction normal to this plane. Specifically, the material matrix remains invariant under rotation by any angle $\theta$ around the axis orthogonal to the plane of isotropy. For clarity of exposition, we assume that the material has symmetric physical properties orthogonal to the $z$ axis, which is normal to the plane of isotropy ($xy$ plane).  In this paper, we characterize the transversely isotropic materials by the following five elastic constants: elastic modulus and shear modulus along the $xy$ plane and along the $z$ axis, and its Poisson ratio along the $xz$ plane. In consequence, we consider the compliance matrix
% \begin{equation} \label{eq:comp_tIsotropic}
%     \mathbf{\bar{S}} = \begin{bmatrix}
%     \frac{1}{E_x} & -\frac{\nu_{x}}{E_x} & -\frac{\nu_{zx}}{E_z} & 0 & 0 & 0 \\
%     -\frac{\nu_{x}}{E_x} & \frac{1}{E_x} & -\frac{\nu_{zx}}{E_z} & 0 & 0 & 0 \\
%     -\frac{\nu_{xz}}{E_x} & -\frac{\nu_{xz}}{E_x} & \frac{1}{E_z} & 0 & 0 & 0 \\
%     0 & 0 & 0 & \frac{1+\nu_x}{E_x} & 0 & 0 \\
%     0 & 0 & 0 & 0 & \frac{1}{2G_{xz}} & 0 \\
%      0 & 0 & 0 & 0 & 0 & \frac{1}{2G_{xz}} 
%     \end{bmatrix},
% \end{equation}
\begin{equation} \label{eq:comp_tIsotropic}
    \mathbf{\bar{S}} = \begin{bmatrix}
 \frac{1}{E_x} & \frac{1}{E_x}-\frac{1}{2G_{xy}} & -\frac{\nu_{xz}}{E_z}  & 0 & 0 & 0 \\
    \frac{1}{E_x}-\frac{1}{2G_{xy}} & \frac{1}{E_x} & -\frac{\nu_{xz}}{E_z}  & 0 & 0 & 0 \\
    -\frac{\nu_{zx}}{E_x} & -\frac{\nu_{zx}}{E_x} & \frac{1}{E_z} & 0 & 0 & 0 \\
    0 & 0 & 0 & \frac{1}{2G_{xz}} & 0 & 0 \\
    0 & 0 & 0 & 0 & \frac{1}{2G_{xz}} & 0 \\
     0 & 0 & 0 & 0 & 0 & \frac{1}{2G_{xy}} 
    \end{bmatrix},
\end{equation}
where $E_x, G_{xy}$ are the Young's modulus and the shear modulus in the $xy$ plane, and $E_z, \nu_{xz},$  and $G_{xz}$ are the  Young's modulus, Poisson ratio, and shear modulus, respectively, in the $z$ direction. Together with the relation $\frac{\nu_{zx}}{E_x} = \frac{\nu_{xz}}{E_z}$, we observe that five independent parameters are needed to define the stress-strain relation for transversely isotropic materials and are denoted in what follows by $$\mathbf{p} = [E_x, E_z, G_{xy}, G_{xz}, \nu_{xz}] \in \R^5.$$ We remark that the elastic constants $E_x, E_z, G_{xy}, G_{xz}$ should have positive values since they represent Young's and shear moduli. On the other hand, the Poisson ratio may be negative. A more explicit formulation of the stress-strain relations for transversely isotropic materials can be obtained by computing the stiffness matrix corresponding to \eqref{eq:comp_tIsotropic}.

\subsection{Solving the eigenvalue problem} \label{sec:solve_COMSOL}
% ---------------------------------------------------------------------------
% 
We use COMSOL Multiphysics\textsuperscript{\textregistered} with its Structural Mechanics Module to obtain a set of eigenfrequencies from a set of parameters $\mathbf{p}$ as input. In other words, COMSOL defines a function $f\colon \mathbf{p} \mapsto (\Lambda_i)_{i=1}^N, \mathds R^5 \to \mathds R^N$, mapping the parameters to the $N$ smallest eigenvalues arranged in ascending order.

In order to set up the computational study in COMSOL, we split the domain $\Omega$ for the engine's rotor into the three components, illustrated in Figure \ref{fig.rotor_components}. We then treat each component as a linear elastic material, using the material properties of structural steel for the shafts and plates and the material properties of copper for the rings and bars. Finally, we assign the transversely isotropic material properties specified by the parameter $\mathbf{p}$ to the rotor core. 

Afterward, we run an eigenvalue study in COMSOL Multiphysics and specify the number of eigenvalues we want to compute. This returns the eigenvalues and their corresponding eigenfunctions. The first six computed eigenvalues correspond to rigid motions and do not provide any information about the parameters. Therefore, we require more than six eigenvalues to obtain eigenvalues corresponding to bending or torsional modes. For this paper, we consider at most three bending pairs and one torsional mode, as this is typically the data provided by most measuring devices. Thus, we require COMSOL to compute the first 13 eigenvalues, including the six eigenvalues corresponding to rigid motions.

\section{Numerical algorithms for parameter recovery}  \label{sec:num_methods}
% ---------------------------------------------------------------------------
% 
In this section, we discuss the inverse problem for recovering the parameters $\mathbf{p} = [E_x, E_z, G_x, G_{xz}, \nu_{xz}]$ from the measured eigenvalues $\Lambda_1,\dots, \Lambda_N$ of \eqref{eq:eigenvalue_problem}. 

% ---------------------------------------------------------------------------
\subsection{Least squares minimization}\label{sec:lsq}
% ---------------------------------------------------------------------------
% 
% In this section, we state the inverse eigenvalue problem for the material parameters $\mathbf{p} \in \R ^5$, and propose some methods for recovering the said parameters. Suppose that a measuring procedure recovers a set of eigenvalues $\Lambda_1,\dots,\Lambda_N$ for the eigenvalue problem \eqref{eq:eigenvalue_problem}. 
% % Suppose further that data about the values of the corresponding eigenfunctions $\varphi_1,\dots \varphi_N$ at certain points $\x_1,\dots, \x_k$ of the domain are also available. 
% Our aim is then to recover the parameters $\mathbf{p} = [E_p, E_z, G_p, G_{zp}, \nu_{pz}]$ from these data. 

% We denote by $f:\R^5\rightarrow\R^N$ the process described in Section \ref{sec:solve_COMSOL} which solves for the first $N$ natural eigenvalues $\Lambda_1,\dots, \Lambda_N$ of \eqref{eq:eigenvalue_problem} corresponding to $\mathbf{p}$. 
% and $g:\R^5 \rightarrow \R^{k\times N}$ the mapping from $\mathbf{p}$ to the (normalized) values of the eigenfunctions $\varphi_1,\dots, \varphi_N$ corresponding to the eigenvalues $\Lambda_1, \dots, \Lambda_N$ at $\x_1,\dots,\x_k$. Of course, $k$ depends on some practical matters, such as the number of sensors available, and the number of measurements that will be made. 
We explore two least squares minimization approaches for the recovery task, depending on whether noise is present in the measurement data. We assert, without proof, that the corresponding finite-dimensional optimization problems have unique minimizers in a neighborhood of the ground truth $\mathbf{p}_{\mathrm{true}}$ and that the iterative algorithms introduced below can be initialized within this neighborhood. Later, our numerical simulations seem to validate this assumption.

% ---------------------------------------------------------------------------
\subsubsection{Cost functions}
% ---------------------------------------------------------------------------
% 
For the noise free measurement data, we will utilize simple least squares minimizer, i.e., we seek $\mathbf{p}$ that solves
\begin{equation}\label{eq:minimisation_exact}
 \min_{\mathbf{p} \in \R^5} \norm{f(\mathbf{p})-\bm{\Lambda}}{2}^2,
\end{equation}
where $\bm{\Lambda} = (\Lambda_i)_{i=1}^N$.
In practise, any measurement contains noise. Here, we consider the case for which we assume that the measurements come with multiplicative random noise. That is, each entry of the measurement $\bm\Lambda^\delta$ is given by $\Lambda_i ^\delta = (1+\zeta_i) \Lambda_i$, where $\zeta_i$ is drawn from a normal distribution $N(0, \delta^2)$ with mean $0$ and variance $\delta^2$, which results in the likelihood distribution $\Lambda_i^\delta | \Lambda_i  \sim N(\Lambda_i, \Lambda_i ^ 2 \delta ^2)$. The corresponding maximum likelihood estimator minimizes
\begin{equation}\label{eq:minimisation_approx}
 \min_{\mathbf{p} \in \R^5} \sum_{i=1}^N \frac{(f(\mathbf{p})_i - \Lambda_i^\delta)^2}{f(\mathbf{p})_i^2 \delta ^2} + \log(f(\mathbf{p})_i).
\end{equation}
Notice that perfect recovery is not possible in this case, and the reconstruction error will depend on the noise level $\delta$.

% ---------------------------------------------------------------------------
\subsubsection{Numerical implementation}
% ---------------------------------------------------------------------------
% 
The optimization problems \eqref{eq:minimisation_exact} and \eqref{eq:minimisation_approx} can be solved via iterative methods. In particular, we employ the minimizers from Python's \code{SciPy} optimization package. The primary input required by the algorithms is an initial guess for the parameters. In addition to the initial guess, we require a priori bounds for the parameters in order to avoid negative or nonphysical solutions. This translates to reformulating \eqref{eq:minimisation_exact} and \eqref{eq:minimisation_approx} as constrained minimization problems. In this paper, we consider Hessian-based methods, such as the L-BFGS-B \cite{BLNZ95-LBFGS,ZBLN97-LBFGSB} or the sequential least squares (SLSQP) \cite{K88-SQP}, which offer a fast convergence for solving the minimization problem. However, some parameters, e.g. the Poisson ratio, have a very small effect on the eigenvalues, and the optimization landscape of \eqref{eq:minimisation_exact} and \eqref{eq:minimisation_approx} is very flat. In this case, we use trust-region methods \cite{CGT00-trust}, which are more robust but come with additional computational costs, for solving the minimization problem. In the next section, we also discuss ensemble Kalman inversion, an alternative derivative-free method for recovering the parameters.

% ---------------------------------------------------------------------------
\subsection{Ensemble Kalman inversion}\label{sec:EKI}
% ---------------------------------------------------------------------------
% 
We now discuss ensemble Kalman inversion (EKI)  as an alternative numerical method for solving the inverse eigenvalue problem associated with \eqref{eq:eigenvalue_problem}. The use of EKI for inverse problems was inspired by the ensemble Kalman filter \cite{E94-enKf}, an ensemble-based algorithm originally designed for high dimensional data assimilation problems. Its main advantage is that it does not require calculating derivatives, which can benefit problems where derivatives are difficult or costly. 

The EKI starts with an initial ensemble of candidate solutions $\{\mathbf{p}_0^{(j)}\}_{j=1}^J \subset \mathbb{R}^5$. At each step, the evolution of the ensemble is determined via a regularization parameter, the sample mean, 
sample covariance, and a perturbation of the measured data. We present in Algorithm \ref{algo:ensemble_up} an adaptation of \cite[Algorithm 1]{I21-adaptive_reg} for the inverse eigenvalue problem considered in this work.

\begin{algorithm}[ht]
	\caption{Ensemble update.}\label{algo:ensemble_up}
	\begin{algorithmic}[1]
        \STATE Given the $n$th ensemble update $\{\mathbf{p}_n^{(j)}\}$, compute the forward map $G_n^{(j)} = f(\mathbf{p}_n^{(j)})$.
        \STATE Compute the regularization parameter $\alpha_n$ and check for convergence.
        \STATE Compute the ensemble mean
        \begin{equation}\nonumber
            \mathbf{\bar{p}}_n = \frac{1}{J} \sum_{j=1}^J \mathbf{p}_n^{(j)},  \qquad \bar{G}_n = \frac{1}{J} \sum_{j=1}^J G_n^{(j)}.
        \end{equation}
        \STATE Compute the sample covariance
        \begin{equation}\nonumber
        \begin{aligned}
                        C_n^{GG} &= \frac{1}{J-1} \sum_{j=1}^J \left(G_n^{(j)}-\bar{G}_n\right) \left(G_n^{(j)}-\bar{G}_n\right)^T \\
                        C_n^{\mathbf{p}G} &= \frac{1}{J-1} \sum_{j=1}^J \left(\mathbf{p}_n^{(j)}-\mathbf{\bar{p}}_n\right) \left(G_n^{(j)}-\bar{G}_n\right)^T
        \end{aligned}
        \end{equation}
        \STATE Update the ensemble 
        \begin{equation}\nonumber
        \begin{aligned}
            \mathbf{p}_{n+1}^{(j)} &= \mathbf{p}_n^{(j)} + C_n^{\mathbf{p}G}\left( C_n^{GG} + \alpha_n\Gamma \right)^{-1}\left(\bm\Lambda^\delta +\sqrt{\alpha_n}\eta_n^{(j)} - G_n^{(j)} \right)
        \end{aligned}
        \end{equation}
        where $\Gamma$ is the covariance matrix and the perturbations $\eta_n^{(j)}$ are i.i.d. collections of vectors indexed by $j,n$ from the normal distribution $N \left(0, \Gamma\right)$.
	\end{algorithmic}
\end{algorithm}
\noindent For this particular multiplicative noise model, we have that $\Gamma$ is a diagonal matrix with entries $\bm\Lambda_i^2\delta^2$. We remark that the ensemble mean $\mathbf{\bar{p}}_{n+1}$ can be informally interpreted as an iterative Tikhonov regularization
applied to a linearization of $f$ through the ensemble. Indeed, it can be shown that the ensemble mean $\mathbf{\bar{p}}_{n+1}$ approximates
the iterative Tikhonov estimator
\begin{equation}\label{eq:EKI_mean}
\begin{aligned}
m_{n+1} &= \frac{1}{2}\mathrm{argmin}_{\mathbf{p}\in S_0}\left\{ \norm{\Gamma^{-1/2}(\bm\Lambda^\delta-f_n-Df_n(\mathbf{p}-m_n))}{}^2 + \alpha_n\norm{\mathcal{C}_n^{-1/2}(\mathbf{p}-m_n) }{}^2 \right\}\\
     &= m_n + \mathcal{C}_nDf_n^*\left(Df_n\mathcal{C}_nDf_n^*+\alpha_n\Gamma \right)^{-1}(\bm\Lambda^\delta-f_n),
     \end{aligned}
     \end{equation}
     with
     \begin{equation}\nonumber
    \mathcal{C}_{n+1} = \mathcal{C}_n - \mathcal{C}_nDf_n^*\left(Df_n\mathcal{C}_nDf_n^*+\alpha_n\Gamma\right)^{-1}Df_n\mathcal{C}_n,
\end{equation}
where $S_0 := \mathrm{span}\{\mathbf{p}_0^{(j)}\}_{j=1}^J$ is the subspace generated by the initial ensemble, $Df$ denotes the Fr\'echet derivative of $f$, $f_n = f(m_n), Df_n = Df(m_n),$ and $Df_n^*$ is the adjoint of $Df$ evaluated at $m_n$. For more detailed discussion, see \cite{I21-adaptive_reg} and references therein. 

In this paper, we take $\alpha_n=1$ for all iterates, which corresponds to the most straightforward implementation of an EKI method. For other variants of EKI, including dynamic choices for $\alpha_n$ and covariance inflation, we refer the readers to \cite{CT19-EKI_conv_acceleration,I21-adaptive_reg}.  

We now discuss the termination criteria for the EKI iterations, starting with the case where there is no noise. In this case, we employ stopping criteria based on the change in the ensemble between successive updates and the updated ensemble's variance. The change in the ensemble is computed via the relative change element-wise, denoted by $\mathbf{r}_{n+1}$, where 
\begin{equation}
    \mathbf{r}_{n+1}^{(j)}: = \frac{|\mathbf{p}_{n+1}^{(j)} - \mathbf{p}_{n}^{(j)}|}{\mathbf{p}_n^{(j)}}.
\end{equation}
Given tolerance values $\mathrm{tol}_c$ and $\mathrm{tol}_v$ for the change and the variance in the ensemble, respectively, we then stop the EKI iterates once $\norm{\mathbf{r}_{n}}{}\leq \mathrm{tol}_c$ and $\mathrm{var}(\mathbf{p}_{n} / \mathbf{\bar{p}}_{n}) \leq \mathrm{tol}_v$. Essentially, the stopping criterion based on $\mathbf{r}$ determines that the ensemble has converged once the change between successive iterates of the ensemble becomes small enough. On the other hand, the stopping criterion based on the variance ensures that upon convergence, the members of the final ensemble should be very close to each other.

When noise is present, we terminate the EKI iteration once the variance of the updated ensemble is less than the measurement variance, i.e., the points of the ensemble have gathered together in one region. This termination procedure is similar to the discrepancy principle \cite{I13-EKI_IP,KNS08-iter}, which terminates the EKI method for the first $n$ such that 
\begin{equation}\label{eq:disc_principle}
\norm{\Gamma^{-1/2}(\bm\Lambda^\delta - f(\mathbf{\bar{p}}_n))}{} \leq \norm{\Gamma^{-1/2}\bm\zeta}{}.
\end{equation}
Here, we see that the termination criterion is natural in that we stop the iterations once the forward map of the ensemble mean is very close to the observed data relative to the noise.

\section{Numerical results} \label{sec:Num_res}
\begin{figure}[t!]
\centering
	\begin{tabular}{cc}
		\includegraphics[width=0.45\linewidth]{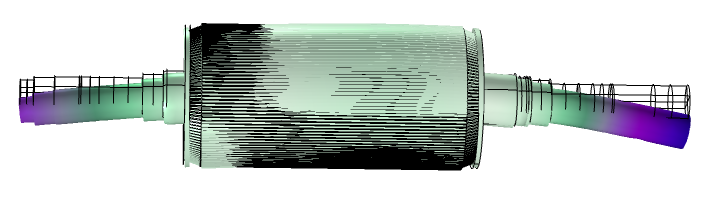} &
		\includegraphics[width=0.45\linewidth]{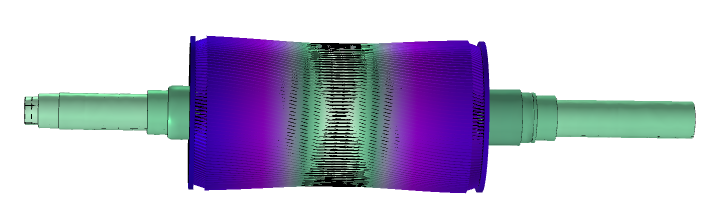} 
	\end{tabular}
	\caption{Illustration of mode types. Left: bending, Right: torsional.} \label{fig.rotor_modes}
\end{figure}

% ---------------------------------------------------------------------------
% 
We start by performing numerical tests to assess how many eigenvalues are needed in order to recover reasonable values for the parameters. For the model considered in this paper, there are two main types of eigenvalues: those that correspond to bending modes and those corresponding to torsional modes. Bending modes usually affect the shafts (see Figure \ref{fig.rotor_modes}, left), whereas torsional modes mostly affect the central region where the rotor core and the rings and bars are located (see Figure \ref{fig.rotor_modes}, right).

For this test, we consider the true parameters $\mathbf{p}_{\mathrm{true}} = [$2e11, 2e8, 7.6923e10, 5e8, 0.3$]$, where the units of the Young's and shear moduli are given in Pascal. We start by considering a scenario where only one parameter is unknown and proceed by considering cases with more unknown parameters. Before performing computations, we remark that having proper bounds for the parameters is important. For most materials, the compliance tensor should be positive semidefinite. Constraints for orthotropic materials that guarantee the positive semidefiniteness of the compliance tensor were provided in \cite{AF87-constraints}. For transversely isotropic materials, these constraints may be rewritten as
% \begin{equation}\label{eq:cons}
% \begin{cases}
%     E_x-\nu_{xz}^2E_z &\geq 0,\\
%     E_x-4G_{xy} &\leq 0, \\
%     E_x-E_x\left| \frac{E_x}{2G_{xy}} -1\right| -2\nu_{xz}^2 E_z & \geq 0
% \end{cases}
% \end{equation}
% 
\begin{equation}\label{eq:cons_simp} %simplified version of the constraints
\begin{cases}
    E_x-\nu_{xz}^2E_z &\geq E_x\left| \frac{E_x}{2G_{xy}} -1\right| + \nu_{xz}^2 E_z,\\
    E_x-4G_{xy} &\leq 0. \\
\end{cases}
\end{equation}

\begin{figure}[t!]\centering
 \includegraphics[width=0.47\linewidth]{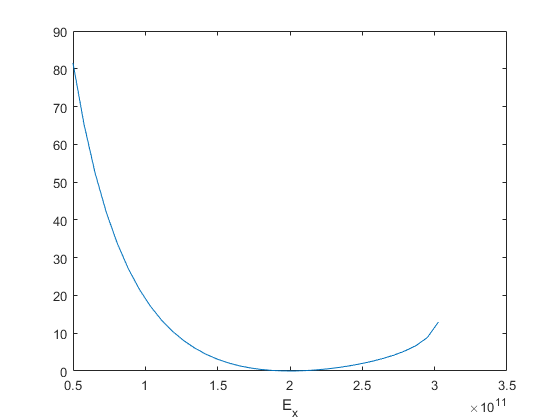} \hfill
 \includegraphics[width=0.47\linewidth]{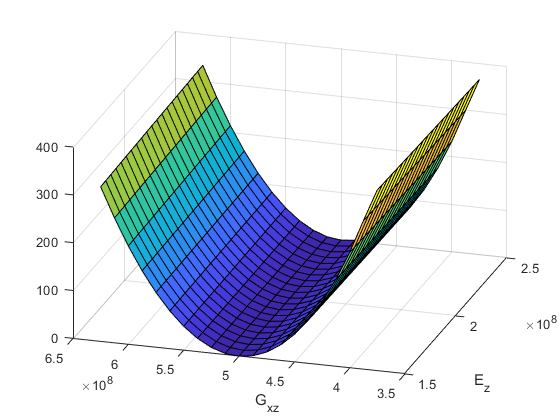}
 \caption{Illustrations of the optimization landscapes given three bending pairs and one torsional mode: The ordinates illustrate the values of the cost function \eqref{eq:minimisation_exact}. In the left picture, the abscissa illustrates the physically reasonable values for $E_x \in [5\text{e}10,3\text{e}11]$. In the right picture, the abscissas cover the physically reasonable values of $E_z$ and $G_{xz}$. All other parameters are assumed to be known exactly.}\label{FIG:landscapes}
\end{figure}

Considering the case where the main unknown is $E_x$, we illustrate the importance of the bounds \eqref{eq:cons_simp}. For the particular case of $\mathbf{p}_{\mathrm{true}}$, it can be computed that in order to satisfy \eqref{eq:cons_simp}, we should have $2.36$e9 $\leq E_x \leq 3.076$e11.  Indeed, choosing parameter values outside the range of \eqref{eq:cons_simp} may yield nonphysical results, e.g. choosing $E_x=3.1$e11,  
%and $E_x=3.33$e11. 
we observed that a purely imaginary eigenvalue was obtained in addition to the six rigid body modes. As there is now a nonphysical mode obtained before the bending and torsional modes, the cost function \eqref{eq:minimisation_exact} will exhibit artificial maxima for this value of $E_x$.
% \AR{This is illustrated in Figure \ref{fig.Ex_cost}, right.}{}{This should be removed}
On the other hand, there is a clear minimum at $E_x=2$e11 in the interval (1e11, 3e11), where the mathematical constraints are satisfied; see the left panel of Figure \ref{FIG:landscapes}. Additionally, we may also use the material properties to provide upper and lower bounds for the parameters; for example, the Young's modulus for steel is expected to be between 190-210GPa \cite{B09-constitutive_eq,AJ12-engmat}, and thus we expect $E_x$ to be within this range of values.
% }{} \AR{it can also be seen in Figure \ref{fig.Ex_cost}, left, that over the interval (1e11, 3e11), as the constraints \eqref{eq:cons_simp} are satisfied, there is a clear global minimum around $E_x=2$e11.}{}{} 

Let us briefly consider the identification of two unknown parameters: $E_z$ and $G_{xz}$. The right panel of Figure \ref{FIG:landscapes} suggests that it is much easier to identify $G_{xz}$ than it is to identify $E_z$ as the graph is (visually) almost constant with respect to $E_z$. This indicates a weak dependence between the cost function value and the parameter $E_z$.
% }{I inserted Figure \ref{fig.rotor_withsens}. Figures \ref{fig.Ex_cost} and \ref{fig.Ez_Gxz} can be removed.}
Similar observations can be made for other combinations of parameter pairs.

% ---------------------------------------------------------------------------
\subsection{Least squares method}\label{sec:lsq_nonoise}
% ---------------------------------------------------------------------------
% 
We present numerical results from solving the optimization problem \eqref{eq:minimisation_exact} via the SLSQP method in Python's \code{SciPy} optimization package. Results obtained from L-BFGS-B and trust constraint methods are very similar and thus have not been reported. If four of the five elastic constants are known, then the remaining elastic constant can be recovered with negligible errors even with only two pairs of bending modes, see Table \ref{tab:noNoise1}. Increasing the amount of available data to three pairs of bending modes improves the estimate's accuracy, whereas adding information on one torsional mode does not significantly impact the accuracy of recovering a single parameter.  

% \begin{table}[h]
%     \centering
% \begin{tabular}{|l|*{3}{c|}}\hline
% \backslashbox{Parameter}{Data}
% &2 bending pairs& 3 bending pairs & \multicolumn{1}{c|}{ 3 bending pairs,}
% \\
% &&& 1 torsional\\\hline
% $E_x$ &2.3513e-4&2.1383e-5& 2.4021e-5\\\hline
% $E_z$ &2.0454e-4&1.1508e-5& 7.0982e-5\\\hline
% $G_{xy}$ &9.0969e-5&4.0289e-4& 1.1535e-4\\\hline
% $G_{xz}$ &7.7437e-6&2.3276e-6& 2.6595e-7\\\hline
% $\nu_{xz}$ &5.3435e-3&7.3308e-3& 5.3729e-3\\\hline
% \end{tabular}
%     \caption{Relative errors in recovering a single parameter from different sets of available data via SLSQP, no noise}
%     \label{tab:noNoise1}
% \end{table}

Next, we explore the recovery of parameter pairs. We assume that three of the five elastic constants are known, and we need to recover two parameters simultaneously. Using only two bending pairs, most parameter pairs can be recovered with an error less than $1\%$, see Table \ref{tab:noNoise2}. Note, however, that for the pair $(E_z,G_{xz})$, we obtain an error of around $15\%$ for $E_z$. Also, pairs that include the Poisson ratio yield errors between $1\%$ and $10\%$. Given information about one more bending pair, the parameter pair $(E_z,G_{xz})$ can now be recovered quite accurately. However, the recovered parameter pair for $(E_x,G_{xy})$ still exhibits errors between $1\%$ and $10\%$. Upon adding information about the first torsional mode, the parameter pair $(E_x,G_{xy})$ can now be recovered with an error less than $1\%$. However, we note that recovered parameter pairs, which include the Poisson ratio, maintain errors between $1\%$ and $10\%$.  One of the reasons for relatively large errors is the failure of the optimizer to converge, e.g. due to the solver encountering incompatible constraints after a certain number of iterates. The presence of more data helps in mitigating this issue. For this particular test, we found that for parameter pairs that include the Poisson ratio, at least five pairs of bending modes are needed for the SLSQP to recover the parameters accurately. %\AR{}{}{In the introduction and Figure \ref{fig.rotor_modes}, we emphasize that the rotor core is (almost) exclusively affected by rotational movements. Thus, it appears natural that we cannot recover its parameters accurately by just considering bending modes. Am I missing anything here?}

% \begin{table}[h]
%     \centering
% \begin{tabular}{|l|*{3}{c|}}\hline
% \backslashbox{Parameter}{Data}
% &2 bending pairs& 3 bending pairs & \multicolumn{1}{c|}{ 3 bending pairs,}
% \\
% &&& 1 torsional\\\hline
% $(E_x,E_z)$ &(5.751e-4, 1.410e-3)&(4.757e-4, 1.005e-3)&(1.045e-5, 1.784e-4)\\\hline
% $(E_x,G_{xy})$ &(1.666e-2, 4.972e-2)&(2.622e-2, 7.143e-2)&(4.913e-5, 7.862e-4)\\\hline
% $(E_x,G_{xz})$ &(6.908e-4, 2.173e-5)&(2.028e-5, 8.952e-9)&(2.213e-6, 6.434e-6)\\\hline
% $(E_x,\nu_{xz})$ &(2.430e-4, 9.635e-2)&(4.924e-5, 2.099e-2)& (1.546e-6, 8.338e-2) \\\hline
% $(E_z,G_{xy})$ &(2.313e-4, 6.270e-4)& (5.709e-4, 4.281e-4)& (3.263e-4, 5.622e-4)\\\hline
% $(E_z,G_{xz})$ &(1.520e-1, 1.748e-3) & (4.753e-3, 3.292e-5)& (4.997e-4, 3.321e-6)\\\hline
% $(E_z,\nu_{xz})$ &(1.665e-4, 2.402e-2) & (1.970e-3, 2.685e-3)& (1.078e-4, 2.135e-2)\\\hline
% $(G_{xy},G_{xz})$ &(4.392e-4, 1.039e-5)& (1.232e-4, 3.708e-6) & (5.929e-4, 6.034e-6)\\\hline
% $(G_{xy},\nu_{xz})$ &(5.782e-4, 1.939e-2)&(6.911e-5, 1.360e-2)& (1.550e-4, 3.242e-2)\\\hline
% $(G_{xz},\nu_{xz})$ &(6.912e-6, 8.238e-2) &(3.410e-6, 8.177e-2)& (3.049e-6, 8.246e-2)\\\hline
% \end{tabular}
% \caption{Relative errors in recovering two parameters from different sets of available data via SLSQP, no noise}
%     \label{tab:noNoise2}
% \end{table}

Upon considering the recovery of parameter triples, we observed that the SLSQP generally fails to converge with only two pairs of bending modes. On the other hand, for triplets excluding the Poisson ratio, the SLSQP converges with three bending pairs and one torsional mode, with errors of order 1e-4. 
% the errors in the Poisson ratio were still huge, presumably due to the fact that the Poisson ratio has a very small effect on the eigenvalues. 
% For instance, the parameter triple $(E_x,G_{xy},G_{xz})$ was estimated with relative error (7.925e-3, 1.892e-1, 1.089e-1), where we see that two out of the three parameters had error larger than $10\%$. Similar observations were noted for other parameter triples recovered via the SLSQP method. 

% ---------------------------------------------------------------------------
\subsection{Ensemble Kalman inversion}\label{sec:EKI_nonoise}
% ---------------------------------------------------------------------------
% 
This section employs the EKI method for solving the inverse eigenvalue problem. Here, we use the EKI method with $J=60$ ensemble members, and the initial ensemble is drawn from a uniform distribution over the interval $(0.5\mathbf{p}_{\mathrm{true}},1.5\mathbf{p}_{\mathrm{true}})$.

We notice that one of the advantages of the EKI method is that it obtains a much better estimate of the parameters compared to the SLSQP method (compare Tables \ref{tab:EKI_noNoise1} and \ref{tab:EKI_noNoise2} with Tables \ref{tab:noNoise1} and \ref{tab:noNoise2}). This is because the default convergence criterion of SLSQP is more relaxed than the one imposed for EKI. The SLSQP can achieve a similar magnitude of accuracy if a stricter convergence criterion is imposed. However, this comes with additional computational cost and takes much longer to converge than the EKI. We now note that the accuracy of EKI comes with a higher computational cost associated with running $J\times N$ forward maps, where $N$ is the number of iterates. This is not a significant issue since the ensemble updates are based on forward evaluations of independent ensemble members, which can be computed in parallel. 

Another interesting observation is that with only 2 bending pairs, the EKI method can recover parameter pairs quite accurately; see Tables \ref{tab:EKI_noNoise1} and \ref{tab:EKI_noNoise2}. In particular, at Table \ref{tab:EKI_noNoise2}, we notice that the only pair which illustrated a vast difference in rate of recovery between using two and three bending pairs was for $(E_x, G_{xz})$. Upon proceeding to explore the accuracy of the EKI method in recovering a more extensive set of parameters, we find that even 2 bending pairs are good enough to recover parameter quadruples. We note, however, that for recovering 3 or more parameters, having more data offers an advantage, in the sense that 2 bending pairs can recover parameters with the error of order 1e-4, whereas 3 bending pairs and one torsional mode can recover parameters with the error of order 1e-6. Finally, in attempting to recover five parameters, we find that 2 bending pairs result in an error of approximately 8.28\% in the recovery of the Poisson ratio and slightly less than 1\% for the other parameters. On the other hand, 3 bending pairs and 1 torsional mode can recover the Poisson ratio with an error of only 0.0026\% and an error less than 0.0001\% for all other parameters. This illustrates that even in the absence of noise, 2 bending pairs are not enough to recover all five parameters simultaneously. This is expected since four data points are generally insufficient to recover five unknowns. 

% ---------------------------------------------------------------------------
\subsection{Noisy data}
% ---------------------------------------------------------------------------
% 
Section \ref{sec:EKI_nonoise} illustrated that the EKI method with $J=60$ ensemble members performs better than the SLSQP for noiseless data. Thus, this section presents numerical tests for the EKI method on data with multiplicative noise. The purpose of the numerical tests is twofold: Firstly, we show that the availability of more data makes the parameters recovered less sensitive to noise. In this aspect, we compare the parameters recovered when only 2 bending pairs are available to the parameters recovered when 3 bending pairs and 1 torsional mode are available. Secondly, we assess the magnitude of allowable noise for which we still recover sensible data by considering increasing noise levels, starting with $\delta = 1$e-10. We illustrate these results with the parameter pairs $(E_x,G_{xy})$ and $(E_z,G_{xz})$. The choice of these parameter pairs is due to the fact that they are the most delicate to recover in the sense that these pairs can only be recovered accurately via the SLSQP when we have at least three bending pairs together with one torsional mode. 

\begin{figure}[ht]
	\begin{tabular}{cc}
\begin{tikzpicture}[scale=0.75]
\begin{loglogaxis}[
xlabel=$\delta$,
ylabel=$E_r$,
legend pos=south east,
 xtick={1e-2,  1e-4, 1e-6, 1e-8, 1e-10},
    xticklabels={$10^{-2}$,  $10^{-4}$, $10^{-6}$, $10^{-8}$, $10^{-10}$}
]

\addplot[mark=o, color=blue] plot coordinates {
    (1e-02 , 1.5587e-01)
    (1e-03 , 6.2242e-02)
    (1e-04 , 2.0038e-02)
 %   (1e-05 , 7.5596e-05)
    (1e-06 , 5.0404e-04)
    (1e-08 , 3.5529e-06)
   (1e-10 , 1.81978e-08)
%    (1e-12 , 4.913e-05)
};

\addplot[mark=o, color=red] plot coordinates {
(1e-02, 1.6022e-01)
(1e-03, 4.0690e-02)
(1e-04, 5.3842e-02)
%(1e-05, 6.7782e-04)
(1e-06, 1.4566e-03)
(1e-08, 1.0450e-05)
(1e-10, 2.6535e-08)
%(1e-12, 7.862e-04)
};

\legend{$E_x$ , $G_{xy}$}
\end{loglogaxis}
\end{tikzpicture}
&
\begin{tikzpicture}[scale=0.75]
\begin{loglogaxis}[
xlabel=$\delta$,
ylabel=$E_r$,
legend pos=south east,
 xtick={1e-2,  1e-4, 1e-6, 1e-8, 1e-10},
    xticklabels={$10^{-2}$,  $10^{-4}$, $10^{-6}$, $10^{-8}$, $10^{-10}$}
]

\addplot[mark=o, color=blue] plot coordinates {
    (1e-02 , 2.8435e-03)
  %  (5e-03 , 3.7094e-02)
    (1e-03 , 2.2163e-01)
 %   (5e-04 , 4.3053e-02)
    (1e-04 , 2.7092e-02)
  %  (1e-05 , 2.0086e-03)
   (1e-06 , 1.9583e-04)
   (1e-08, 1.6077e-05)
   (1e-10, 2.7913e-07)
%    (1e-12 , 7.0323e-09)
%	(5.525e-03,  1.121e-05)
};

\addplot[mark=o, color=red] plot coordinates {
(1e-02, 1.5584e-03)
%(5e-03, 4.0768e-04)
(1e-03, 2.2155e-03)
%(5e-04, 3.5969e-04)
(1e-04, 3.0796e-04)
%(1e-05, 2.8786e-05)
(1e-06, 2.0355e-06)
(1e-08, 2.1393e-07)
(1e-10, 2.8937e-09)
%(1e-12, 1.109e-10) 
};

$\legend{$E_z$, $G_{xz}$}$
\end{loglogaxis}
\end{tikzpicture}
\end{tabular}
	\caption{Relative errors for different noise levels, 2 bending pairs. Left: $(E_x,G_{xy})$; right: $(E_z,G_{xz})$. } \label{fig.err_plots1}
 \end{figure}
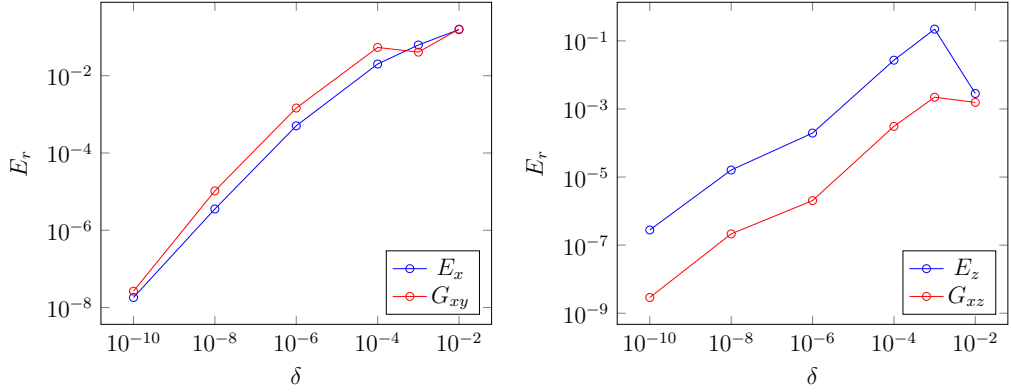

Figure \ref{fig.err_plots1} shows the relative errors in recovering the parameter pairs with different levels of noise $\delta$. These were obtained via the EKI method (see Algorithm \ref{algo:ensemble_up}) with $J=60$ members, under the assumption that only 2 bending pairs are available. The initial ensemble is drawn from a uniform distribution over the interval $(0.5\mathbf{p}_{\mathrm{true}},1.5\mathbf{p}_{\mathrm{true}})$, and the relative errors are measured by comparing the mean of the final ensemble with the true parameter value. That is, if there are a total of $N$ ensemble updates, then the relative errors are computed with
\begin{equation}\nonumber
E_r = \frac{|\mathbf{p}_{\mathrm{true}} - \bar{\mathbf{p}}_N|}{\mathbf{p}_{\mathrm{true}}}.
\end{equation}

%the parameters recovered by the algorithm already have a relative error which is  $10\%$ or larger.
\begin{figure}[ht]
	\begin{tabular}{cc}
\begin{tikzpicture}[scale=0.75]
\begin{loglogaxis}[
xlabel=$\delta$,
ylabel=$E_r$,
legend pos=south east,
 xtick={1e-2,  1e-4, 1e-6, 1e-8, 1e-10},
    xticklabels={$10^{-2}$,  $10^{-4}$, $10^{-6}$, $10^{-8}$, $10^{-10}$}
]

\addplot[mark=o, color=blue] plot coordinates {
    (1e-02 , 5.2179e-02)
    (1e-03 , 2.6849e-02)
    (1e-04 , 1.3799e-03)
 %   (1e-05 , 7.5596e-05)
    (1e-06 , 1.1720e-05)
    (1e-08 , 1.4667e-07)
    (1e-10 , 1.0678e-09)
%    (1e-12 , 4.913e-05)
};

\addplot[mark=o, color=red] plot coordinates {
(1e-02, 2.1039e-01)
(1e-03, 1.3511e-01)
(1e-04, 1.0868e-02)
%(1e-05, 6.7782e-04)
(1e-06, 7.0699e-05)
(1e-08, 1.1859e-06)
(1e-10, 1.3662e-08)
%(1e-12, 7.862e-04)
};

\addplot[mark=o, color=blue, dashed] plot coordinates {
    (1e-02 , 1.5587e-01)
    (1e-03 , 6.2242e-02)
    (1e-04 , 2.0038e-02)
 %   (1e-05 , 7.5596e-05)
    (1e-06 , 5.0404e-04)
    (1e-08 , 3.5529e-06)
  (1e-10 , 1.81978e-08)
%    (1e-12 , 4.913e-05)
};

\addplot[mark=o, color=red, dashed] plot coordinates {
(1e-02, 1.6022e-01)
(1e-03, 4.0690e-02)
(1e-04, 5.3842e-02)
%(1e-05, 6.7782e-04)
(1e-06, 1.4566e-03)
(1e-08, 1.0450e-05)
(1e-10, 2.6535e-08)
%(1e-12, 7.862e-04)
};

\legend{$E_x$ , $G_{xy}$}
\end{loglogaxis}
\end{tikzpicture}
&
\begin{tikzpicture}[scale=0.75]
\begin{loglogaxis}[
xlabel=$\delta$,
ylabel=$E_r$,
legend pos=south east,
 xtick={1e-2,  1e-4, 1e-6, 1e-8, 1e-10},
    xticklabels={$10^{-2}$,  $10^{-4}$, $10^{-6}$, $10^{-8}$, $10^{-10}$}
]

\addplot[mark=o, color=blue] plot coordinates {
    (1e-02 , 5.3008e-02)
  %  (5e-03 , 3.7094e-02)
    (1e-03 , 1.3408e-01)
 %   (5e-04 , 4.3053e-02)
    (1e-04 , 2.6943e-02)
    (1e-05 , 2.0086e-03)
   (1e-06 , 4.7387e-04)
   (1e-08, 8.1277e-06)
   (1e-10, 5.0751e-08)
%    (1e-12 , 4.997e-04)
%	(5.525e-03,  1.121e-05)
};

\addplot[mark=o, color=red] plot coordinates {
(1e-02, 7.9675e-03)
%(5e-03, 4.0768e-04)
(1e-03, 1.7295e-03)
%(5e-04, 3.5969e-04)
(1e-04, 3.9974e-04)
(1e-05, 2.8786e-05)
(1e-06, 3.9099e-06)
(1e-08, 9.5194e-08)
(1e-10, 4.9358e-10)
%(1e-12, 3.321e-06)
};

\addplot[mark=o, color=blue, dashed] plot coordinates {
    (1e-02 , 2.8435e-03)
  %  (5e-03 , 3.7094e-02)
    (1e-03 , 2.2163e-01)
 %   (5e-04 , 4.3053e-02)
    (1e-04 , 2.7092e-02)
  %  (1e-05 , 2.0086e-03)
   (1e-06 , 1.9583e-04)
   (1e-08, 1.6077e-05)
   (1e-10, 2.7913e-07)
%    (1e-12 , 7.0323e-09)
%	(5.525e-03,  1.121e-05)
};

\addplot[mark=o, color=red, dashed] plot coordinates {
(1e-02, 1.5584e-03)
%(5e-03, 4.0768e-04)
(1e-03, 2.2155e-03)
%(5e-04, 3.5969e-04)
(1e-04, 3.0796e-04)
%(1e-05, 2.8786e-05)
(1e-06, 2.0355e-06)
(1e-08, 2.1393e-07)
(1e-10, 2.8937e-09)
%(1e-12, 1.109e-10) 
};

$\legend{$E_z$, $G_{xz}$}$
\end{loglogaxis}
\end{tikzpicture}
\end{tabular}
	\caption{Relative errors for different noise levels, 3 bending pairs, 1 torsional mode. Left: $(E_x,G_{xy})$; right: $(E_z,G_{xz})$. } \label{fig.err_plots2}
 \end{figure}
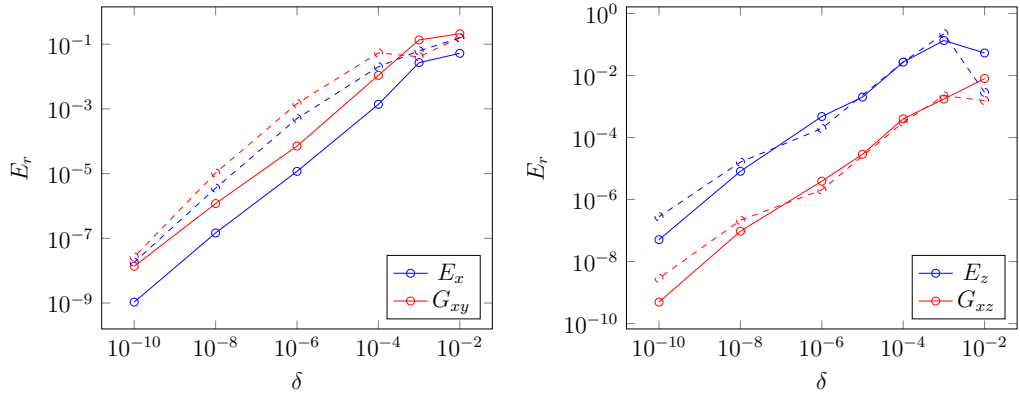

In Figure \ref{fig.err_plots1}, we observe that errors larger than $10\%$ were obtained for noise level $\delta = 1$e-2. We also see that, as expected, the error decreases as the noise level decreases. We now proceed in Figure \ref{fig.err_plots2} to perform a comparison of the results obtained for the same unknowns, but under the assumption that 3 bending pairs and one torsional mode are available. Here, we also observe that lower noise levels result in lower errors in the recovered parameters. Moreover, upon comparing the solid lines with the dashed lines, we observe that the parameters recovered with more data are less noise-sensitive.

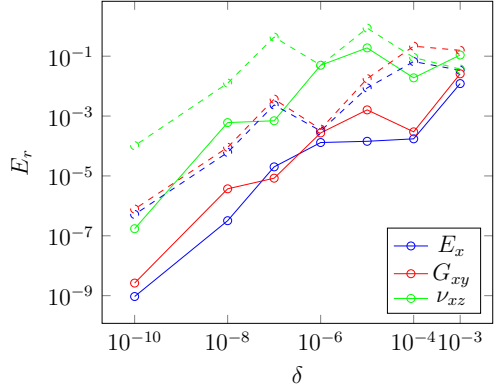
\begin{figure}[ht]
\begin{tikzpicture}[scale=0.75]
\begin{loglogaxis}[
xlabel=$\delta$,
ylabel=$E_r$,
legend pos=south east,
 xtick={1e-3,  1e-4, 1e-6, 1e-8, 1e-10},
    xticklabels={$10^{-3}$,  $10^{-4}$, $10^{-6}$, $10^{-8}$, $10^{-10}$}
]

\addplot[mark=o, color=blue] plot coordinates {
 %   (1e-02 , 5.3008e-02)
  %  (5e-03 , 3.7094e-02)
    (1e-03 , 1.2147e-02)
 %   (5e-04 , 4.3053e-02)
    (1e-04 , 1.7334e-04)
    (1e-05 , 1.4338e-04)
   (1e-06 , 1.3051e-04)
   (1e-07, 2.0060e-05)
   (1e-08, 3.1944e-07)
   (1e-10, 9.3562e-10)
%    (1e-12 , 4.997e-04)
%	(5.525e-03,  1.121e-05)
};

\addplot[mark=o, color=red] plot coordinates {
%(1e-02, 2.5776e-02)
%(5e-03, 4.0768e-04)
(1e-03, 2.5776e-02)
%(5e-04, 3.5969e-04)
(1e-04, 2.9683e-04)
(1e-05, 1.6075e-03)
(1e-06, 2.7315e-04)
(1e-07, 8.3976e-06)
(1e-08, 3.7059e-06)
(1e-10, 2.6301e-09)
%(1e-12, 3.321e-06)
};

\addplot[mark=o, color=green] plot coordinates {
%(1e-02, 7.9675e-03)
%(5e-03, 4.0768e-04)
(1e-03, 1.0940e-01)
%(5e-04, 3.5969e-04)
(1e-04, 1.8896e-02)
(1e-05, 1.8920e-01)
(1e-06, 4.9810e-02)
(1e-07, 6.9544e-04)
(1e-08, 6.0490e-04)
(1e-10, 1.7244e-07)
%(1e-12, 3.321e-06)
};
\addplot[mark=o, color=blue, dashed] plot coordinates {
 %   (1e-02 , 2.8435e-03)
  %  (5e-03 , 3.7094e-02)
    (1e-03 , 3.3852e-02)
 %   (5e-04 , 4.3053e-02)
    (1e-04 , 6.6096e-02)
    (1e-05 , 8.6937e-03)
   (1e-06 , 3.1726e-04)
   (1e-07, 2.3498e-03)
   (1e-08, 6.0868e-05)
   (1e-10, 5.1560e-07)
%   (1e-10, 2.7913e-07)
%    (1e-12 , 7.0323e-09)
%	(5.525e-03,  1.121e-05)
};

\addplot[mark=o, color=red, dashed] plot coordinates {
%(1e-02, 1.5584e-03)
%(5e-03, 4.0768e-04)
(1e-03, 1.5685e-01)
%(5e-04, 3.5969e-04)
(1e-04, 2.1576e-01)
(1e-05, 1.7708e-02)
(1e-06, 3.6495e-04)
(1e-07, 3.6710e-03)
(1e-08, 9.0080e-05)
(1e-10, 7.8110e-07)
%(1e-10, 2.8937e-09)
%(1e-12, 1.109e-10) 
};

\addplot[mark=o, color=green, dashed] plot coordinates {
%(1e-02, 1.5584e-03)
%(5e-03, 4.0768e-04)
(1e-03, 3.7092e-02)
%(5e-04, 3.5969e-04)
(1e-04, 9.0544e-02)
(1e-05, 8.6543e-01)
(1e-06, 4.9900e-02)
(1e-07, 4.3397e-01)
(1e-08, 1.2551e-02)
(1e-10, 1.0208e-04)
%(1e-10, 2.8937e-09)
%(1e-12, 1.109e-10) 
};

$\legend{$E_x$, $G_{xy}$, $\nu_{xz}$}$
\end{loglogaxis}
\end{tikzpicture}
	\caption{Relative errors for different noise levels $(E_x,G_{xy},\nu_{xz})$. Solid lines: three bending pairs, one torsional mode; Dashed: two bending pairs} \label{fig.err_plots3}
 \end{figure}

\begin{figure}[ht]
\begin{tikzpicture}[scale=0.75]
\begin{loglogaxis}[
xlabel=$\delta$,
ylabel=$E_r$,
legend pos=south east,
 xtick={1e-3,  1e-4, 1e-6, 1e-8, 1e-10},
    xticklabels={$10^{-3}$,  $10^{-4}$, $10^{-6}$, $10^{-8}$, $10^{-10}$}
]

\addplot[mark=o, color=blue] plot coordinates {
 %   (1e-02 , 5.3008e-02)
  %  (5e-03 , 3.7094e-02)
    (1e-03 , 2.1716e-02)
 %   (5e-04 , 4.3053e-02)
    (1e-04 , 2.7975e-03)
    (1e-05 , 1.8090e-05)
   (1e-06 , 7.1631e-05)
 %  (1e-07, 2.0060e-05)
   (1e-08, 1.0349e-05)
   (1e-10, 2.0179e-07)
%    (1e-12 , 4.997e-04)
%	(5.525e-03,  1.121e-05)
};

\addplot[mark=o, color=red] plot coordinates {
%(1e-02, 2.5776e-02)
%(5e-03, 4.0768e-04)
(1e-03, 1.5077e-01)
%(5e-04, 3.5969e-04)
(1e-04, 5.3233e-02)
(1e-05, 5.4744e-03)
(1e-06, 1.5584e-03)
%(1e-07, 8.3976e-06)
(1e-08, 4.3175e-04)
(1e-10, 8.3179e-06)
%(1e-12, 3.321e-06)
};

\addplot[mark=o, color=green] plot coordinates {
%(1e-02, 7.9675e-03)
%(5e-03, 4.0768e-04)
(1e-03, 1.0298e-01)
%(5e-04, 3.5969e-04)
(1e-04, 2.5763e-02)
(1e-05, 7.9620e-04)
(1e-06, 6.5032e-04)
%(1e-07, 6.9544e-04)
(1e-08, 1.3283e-04)
(1e-10, 2.5390e-06)
%(1e-12, 3.321e-06)
};
\addplot[mark=o, color=violet] plot coordinates {
 %   (1e-02 , 2.8435e-03)
  %  (5e-03 , 3.7094e-02)
    (1e-03 , 4.4771e-03)
 %   (5e-04 , 4.3053e-02)
    (1e-04 , 4.2213e-05)
    (1e-05 , 2.7217e-05)
   (1e-06 , 8.2146e-06)
   %(1e-07, 2.3498e-03)
   (1e-08, 3.0701e-06)
   (1e-10, 5.9971e-08)
%   (1e-10, 2.7913e-07)
%    (1e-12 , 7.0323e-09)
%	(5.525e-03,  1.121e-05)
};

\addplot[mark=o, color=brown] plot coordinates {
%(1e-02, 1.5584e-03)
%(5e-03, 4.0768e-04)
(1e-03, 1.4695e-01)
%(5e-04, 3.5969e-04)
(1e-04, 1.0503e-01)
(1e-05, 1.3388e-01)
(1e-06, 1.9143e-01)
%(1e-07, 3.6710e-03)
(1e-08, 3.6194e-02)
(1e-10, 6.8386e-04)
%(1e-10, 2.8937e-09)
%(1e-12, 1.109e-10) 
};

$\legend{$E_x$, $E_z$,$G_{xy}$,$G_{xz}$, $\nu_{xz}$}$
\end{loglogaxis}
\end{tikzpicture}
 \caption{Relative errors for different noise levels, simultaneous recovery of all parameters: three bending pairs, one torsional mode.%; Dashed: two bending pairs
 } \label{fig.err_plots4}
\end{figure}
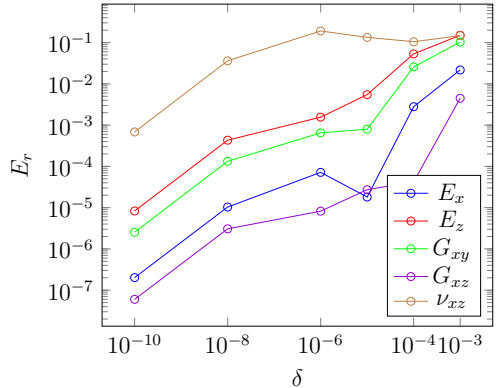

In Figures \ref{fig.err_plots1} and \ref{fig.err_plots2}, the effect of noise on the parameter pairs recovered via 2 bending pairs or 3 bending pairs and 1 torsional mode is not very prominent. Thus, in order to illustrate that the parameter recovery process with 3 bending pairs and 1 torsional mode is less sensitive to noise, we consider the recovery of the parameter triples $(E_x,G_{xy},\nu_{xz})$, $(E_z,G_{xz},\nu_{xz})$. Indeed, it can be seen in Figure \ref{fig.err_plots3} that the errors for the parameters recovered with only 2 bending pairs are much larger than those recovered with more data. Finally, Figure \ref{fig.err_plots4} shows that if we want to recover all five parameters accurately, the noise level needs to be minimal (order 1e-10). On the other hand, taking a noise level of $1$e-6 yields reasonable estimates for four out of the five parameters, with significant errors associated with the estimated Poisson ratio. 

\section{Summary and future work} \label{sec:summary}

In this work, we studied numerical algorithms for parameter recovery in eigenvalue problems for linear elasticity. In particular, we are interested in an application that recovers the elastic constants of a rotor core from measured eigenfrequencies. We assumed that the rotor core is transversely isotropic, i.e., five elastic constants can characterize it. We started by looking at the feasibility of the problem by trying to recover the elastic constants under perfect conditions, i.e., we assume that the numerical model is exact and that no errors were made in the measurement. Here, we made a comparison between the performance of SLSQP and EKI methods. Through this test, we found that in order for the SLSQP to simultaneously estimate parameter pairs accurately (excluding the Poisson ratio), we require at least three pairs of bending modes, together with a torsional mode. On the other hand, with an ensemble size of 60, the EKI can recover up to four parameters simultaneously, even with only two bending pairs. In order to recover all five parameters accurately, we require three bending pairs and one torsional mode.

Following this, we studied a test case with noisy measurements, in which we assumed that the measured eigenfrequencies contain multiplicative noise. Numerical tests for parameter recovery were performed via an EKI method with 60 ensemble members. The results show that for noise levels below $0.01\% $, the parameter pairs can be recovered quite well, with errors close to or much less than $1\%$. Also, as expected, the errors in the recovered parameters are larger for higher noise levels. 

The results reported in this work look promising and opens several natural extensions that may be considered in future work. Firstly, any inversion algorithm requires performing the forward mathematical solver several times. The solution process for computing the eigenvalues is quite expensive and is one of the main bottlenecks of the parameter recovery procedure. Thus, exploring subspace methods for reducing the computational time needed to obtain the eigenvalues will also be interesting. Secondly, we have found that the convergence of ensemble-based inversion methods depends on the choice of the initialization, i.e. constraining the initial ensemble. A natural avenue for future research would be to try reducing the number of iterations needed for convergence by combining Hessian-based and derivative-free methods. Finally, a natural next step would be to demonstrate that the method can be efficient and reliable in the context of industrial applications. This effort would benefit from more realistic mathematical model involving e.g. spatially dependent material parameters.
%Afterward, we look into how the schemes considered in this paper may be adapted and used in the improved model.
% that is as close as possible to the experimental model. The current mathematical model has an accurate depiction of the bending modes but not of the torsional modes. Hence, a natural path to pursue is to look at some ideas that may help capture the torsional mode, such as the usage of an orthotropic material instead of a transversely isotropic material in the model or the assumption that the elastic parameters are non-constant functions, among others. 

\bibliographystyle{alphaAR}
\bibliography{abb_lut}

\appendix
\newpage
\section{Tables with simulation results}

\rowcolors{2}{gray!25}{white}
\begin{table}[ht!]\centering
\begin{tabular}{l|*{3}c}
\toprule
\backslashbox[3cm]{Param.}{Data}
&2 bending pairs& 3 bending pairs & \begin{minipage}{3cm} \centering 3 bending pairs \\ 1 torsional\end{minipage}\\
\midrule
$E_x$ &2.3513e-4&2.1383e-5& 2.4021e-5\\
$E_z$ &2.0454e-4&1.1508e-5& 7.0982e-5\\
$G_{xy}$ &9.0969e-5&4.0289e-4& 1.1535e-4\\
$G_{xz}$ &7.7437e-6&2.3276e-6& 2.6595e-7\\
$\nu_{xz}$ &5.3435e-3&7.3308e-3& 5.3729e-3\\
\bottomrule
\end{tabular}
\caption{Relative errors in recovering a single parameter from different sets of available data via SLSQP, no noise}\label{tab:noNoise1}
\end{table}

\rowcolors{2}{gray!25}{white}
\begin{table}[ht!]\centering
\begin{tabular}{l|*{3}c}
\toprule
\backslashbox[3cm]{Param.}{Data}
&2 bending pairs& 3 bending pairs & \begin{minipage}{3cm} \centering 3 bending pairs \\ 1 torsional\end{minipage}\\
\midrule
$(E_x,E_z)$ &(5.751e-4, 1.410e-3)&(4.757e-4, 1.005e-3)&(1.045e-5, 1.784e-4)\\
$(E_x,G_{xy})$ &(1.666e-2, 4.972e-2)&(2.622e-2, 7.143e-2)&(4.913e-5, 7.862e-4)\\
$(E_x,G_{xz})$ &(6.908e-4, 2.173e-5)&(2.028e-5, 8.952e-9)&(2.213e-6, 6.434e-6)\\
$(E_x,\nu_{xz})$ &(2.430e-4, 9.635e-2)&(4.924e-5, 2.099e-2)& (1.546e-6, 8.338e-2) \\
$(E_z,G_{xy})$ &(2.313e-4, 6.270e-4)& (5.709e-4, 4.281e-4)& (3.263e-4, 5.622e-4)\\
$(E_z,G_{xz})$ &(1.520e-1, 1.748e-3) & (4.753e-3, 3.292e-5)& (4.997e-4, 3.321e-6)\\
$(E_z,\nu_{xz})$ &(1.665e-4, 2.402e-2) & (1.970e-3, 2.685e-3)& (1.078e-4, 2.135e-2)\\
$(G_{xy},G_{xz})$ &(4.392e-4, 1.039e-5)& (1.232e-4, 3.708e-6) & (5.929e-4, 6.034e-6)\\
$(G_{xy},\nu_{xz})$ &(5.782e-4, 1.939e-2)&(6.911e-5, 1.360e-2)& (1.550e-4, 3.242e-2)\\
$(G_{xz},\nu_{xz})$ &(6.912e-6, 8.238e-2) &(3.410e-6, 8.177e-2)& (3.049e-6, 8.246e-2)\\
\bottomrule
\end{tabular}
\caption{Relative errors in recovering two parameters from different sets of available data via SLSQP, no noise}
    \label{tab:noNoise2}
\end{table}

\rowcolors{2}{gray!25}{white}
\begin{table}[ht!]\centering
\begin{tabular}{l|*{3}c}
\toprule
\backslashbox[3cm]{Param.}{Data}
&2 bending pairs& 3 bending pairs & \begin{minipage}{3cm} \centering 3 bending pairs \\ 1 torsional\end{minipage}\\
\midrule
$E_x$ &2.5520e-8&2.6035e-6& 7.5987e-10\\
$E_z$ &7.0336e-10&1.5015e-9& 1.0681e-10 \\
$G_{xy}$ &2.2463e-10&1.6488e-9&2.5560e-10 \\
$G_{xz}$ &4.3389e-10&1.6049e-11&4.0796e-12 \\
$\nu_{xz}$ &3.3416e-8&2.8310e-7&1.1334e-7 \\
\bottomrule
\end{tabular}
    \caption{Relative errors in recovering a single parameter from different sets of available data via EKI, no noise}
    \label{tab:EKI_noNoise1}
\end{table}

\rowcolors{2}{gray!25}{white}
\begin{table}[ht!]\centering
\begin{tabular}{l|*{3}c}
\toprule
\backslashbox[3cm]{Param.}{Data}
&2 bending pairs& 3 bending pairs & \begin{minipage}{3cm} \centering 3 bending pairs \\ 1 torsional\end{minipage}\\
\midrule
$(E_x,E_z)$ &(6.975e-10, 2.191e-8)&(2.359e-7, 6.011e-8)&(3.529e-9, 2.810e-8)\\
$(E_x,G_{xy})$ &(8.948e-9, 2.529e-8)&(1.930e-10, 4.178e-10)&(9.227e-12, 2.519e-10)\\
$(E_x,G_{xz})$ &(1.813e-5, 3.574e-6)&(1.682e-10, 1.942e-12)&(7.346e-11, 8.047e-13)\\
$(E_x,\nu_{xz})$ &(8.023e-10, 2.595e-7)&(2.990e-10, 2.007e-7)& (5.851e-11, 1.318e-7) \\
$(E_z,G_{xy})$ &(3.073e-8, 4.404e-7)& (5.505e-10, 5.178e-10)& (1.746e-9, 1.338e-9)\\
$(E_z,G_{xz})$ &(7.032e-9, 1.109e-10) & (1.069e-9, 7.999e-11)& (1.667e-8, 1.050e-10)\\
$(E_z,\nu_{xz})$ &(1.194e-11, 1.204e-7) & (6.317e-11, 2.558e-8)& (1.077e-10, 1.398e-7)\\
$(G_{xy},G_{xz})$ &(5.441e-10, 9.683e-12)& (9.287e-10, 5.119e-12) & (2.037e-7, 1.011e-10)\\
$(G_{xy},\nu_{xz})$ &(8.103e-10, 2.436e-7)&(8.912e-11, 6.173e-8)& (2.697e-9, 1.373e-7)\\
$(G_{xz},\nu_{xz})$ &(5.159e-12, 6.548e-8) &(7.960e-13, 7.786e-8)& (1.051e-12, 1.717e-7)\\
\bottomrule
\end{tabular}
\caption{Relative errors in recovering two parameters from different sets of available data via EKI, no noise}
    \label{tab:EKI_noNoise2}
\end{table}

%\AR{}{}{I changed the visual appearance of the tables. This makes them fit on fewer pages!}

% 
\end{document}